\documentclass{article}

\newtheorem{fed}{\textbf{Definition}}[section]

\newtheorem{lemma}[fed]{\textbf{Lemma}}
\newtheorem{ex}[fed]{\textbf{Example}}
\newtheorem{rem}[fed]{\textbf{Remark}}

\newtheorem{cor}[fed]{\textbf{Corollary}}
\usepackage{amssymb,bbm,graphicx,epsfig,psfrag,epic,eepic,latexsym}
\usepackage{amsmath}
\usepackage{mathrsfs}
\usepackage[dvips]{color}

\newcommand{\N}{\mathbb{N}}

\newcommand{\Z}{\mathbb{Z}}
\newcommand{\R}{\mathbb{R}}
\newcommand{\C}{\mathbb{C}}

\newcommand{\dom}{\mathrm{dom}}

\begin{document}
\title{Fractal scale Hilbert spaces and scale Hessian operators}
\author{Urs Frauenfelder\footnote{Department of Mathematics and
Research Institute of Mathematics, Seoul National University}}
\maketitle
\begin{abstract}
Scale spaces were defined by H.\,Hofer, K.\,Wysocki, and
E.\,Zehnder. In this note we introduce a subclass of scale spaces
and explain why we believe that this subclass is the right class for
a general setup of Floer theory.
\end{abstract}

\tableofcontents

\section[Introduction]{Introduction}
The definition of a scale Hilbert structure is due to H.\,Hofer,
K.\,Wysocki, and E.\,Zehnder
\cite{hofer1,hofer-wysocki-zehnder1,hofer-wysocki-zehnder2}.
\begin{fed}\label{sca}
A \emph{scale Hilbert space} is a tuple
$$\mathcal{H}=\big\{\big(H_k, \langle \cdot, \cdot \rangle_k\big)\big\}_{k
\in \N_0}$$ where for each $k \in \N_0$ the pair $\big(H_k,\langle
\cdot, \cdot \rangle_k\big)$ is a real Hilbert space and the vector
spaces $H_k$ build a nested sequence
 $H=H_0 \supset H_1 \supset H_2 \supset \ldots$
such that the following two axioms hold.
\begin{description}
 \item[(i)] For each $k \in \N$ the inclusion
 $\big(H_k, \langle \cdot, \cdot \rangle_k\big) \hookrightarrow
 \big(H_{k-1}, \langle \cdot, \cdot \rangle_{k-1}\big)$ is compact.
 \item[(ii)] For each $k \in \N_0$ the subspace
 $H_\infty=\bigcap_{n=0}^\infty H_n$ is dense in $H_k$ with respect
 to the topology induced from $\langle \cdot, \cdot \rangle_k$.
\end{description}
\end{fed}
\begin{rem}
Scale structures can as well be defined on Banach spaces but in this
paper we restrict ourself to the case of Hilbert spaces.
\end{rem}
We further recall the notion of scale continuous map from
\cite{hofer1,hofer-wysocki-zehnder1,hofer-wysocki-zehnder2}. Suppose
that $\mathcal{H}=\{(H_k,\langle \cdot, \cdot \rangle_k)\}$ and
$\mathcal{H}'=\{(H_k',\langle \cdot, \cdot \rangle'_k)\}$ are two
scale Hilbert spaces.
\begin{fed}
A map $f \colon H_0 \to H_0'$ is called \emph{scale continuous} if
for each $k$ the map $f$ restricts to a continuous map $f_k \colon
H_k \to H_k'$.
\end{fed}
\begin{rem}
Using the same notion one can also define scale continuity for maps
which are only defined on an open subset of $H_0$.
\end{rem}
There are two equivalence relations for scale Hilbert spaces which
we explain next.
\begin{fed}A linear, scale continuous map $\Phi$ between
scale Hilbert spaces $\mathcal{H}$ and $\mathcal{H}'$ is called a
\emph{scale isomorphism} if for each $k \in \N_0$ its restriction
$\Phi_k \colon H_k \to H_k'$ is a bijection.
\end{fed}
\begin{rem}
It follows from the open mapping theorem that the inverse of a scale
isomorphism is a scale isomorphism as well.
\end{rem}
\begin{fed}
A \emph{scale isometry} $\Phi$ from $\mathcal{H}$ to $\mathcal{H}'$
is a bijective linear map which restricts for each $k$ to an
isometry $\Phi_k$ between the Hilbert spaces $H_k$ and $H_k'$.
\end{fed}
\begin{fed}
Two scale Hilbert spaces $\mathcal{H}$ and $\mathcal{H}'$ are called
\emph{scale isomorphic}, if there exists a scale isomorphism from
$\mathcal{H}$ to $\mathcal{H}'$. They are called \emph{scale
isometric}, if there exists a scale isometry between them.
\end{fed}
\begin{rem}
Since each scale isometry is a special case of a scale isomorphism,
scale isometric scale Hilbert spaces are automatically scale
isomorphic.
\end{rem}
We also recall from \cite{frauenfelder} the notion of a scale
Hilbert pair.
\begin{fed}
A \emph{scale Hilbert pair} is a pair of Hilbert spaces
$\mathcal{H}=(H_0,H_1)$ such that $H_1 \subset H_0$ is a dense
subset and the inclusion $H_1 \hookrightarrow H_0$ is compact.
\end{fed}
The notions of scale isomorphism and scale isometry for scale
Hilbert pairs is the same as the one for scale Hilbert spaces.

An example of a scale Hilbert pair is the following. Let $f \colon
\N \to (0,\infty)$ be a monotone unbounded function. By $\ell^2$ we
refer as usual to the Hilbert space of square summable sequences. We
say that a sequence $x=(x_1,x_2,\ldots)$ is in $\ell^2_f$ if
$$\sum_{\nu=1}^\infty f(\nu)x_\nu^2<\infty.$$
We give $\ell^2_f$ the structure of a Hilbert space by introducing
for $x,y \in \ell^2_f$ the inner product
$$\langle x, y \rangle_f=\sum_{\nu=1}^\infty f(\nu)x_\nu y_\nu.$$
Then the pair $(\ell^2,\ell^2_f)$ is a scale Hilbert pair. It can be
shown, see \cite{frauenfelder}, that each infinite dimensional scale
Hilbert pair is scale isomorphic to a pair $(\ell^2,\ell^2_f)$. In
particular, each Hilbert space $H_k$ arising in an infinite
dimensional scale Hilbert space $\mathcal{H}$ is separable and hence
isometric to $\ell^2$. This reduces the geography problem for scale
Hilbert spaces to the problem of how infinitely many $\ell^2$-spaces
can be nested into each other to produce a scale Hilbert space.
\\
For $f$ as above we introduce the scale Hilbert space $\ell^{2,f}$
given by
$$\ell^{2,f}_k=\ell^2_{f^k}, \quad k \in \N_0.$$
\begin{rem}
The standard orthogonal basis of $\ell^2$ is a common orthonormal
basis of $\ell^{2,f}_k$ for $k \in \N_0$. This is a nontrivial issue
already in finite dimensions. While two scalar products always admit
a common orthogonal basis this in general fails for three scalar
products.
\end{rem}
We are now in position to define the notion of a fractal scale
Hilbert space.
\begin{fed}
A scale Hilbert space $\mathcal{H}$ is called \emph{fractal} if
there exists a unbounded monotone function $f \colon \N \to
(0,\infty)$ such that $\mathcal{H}$ is scale isomorphic to
$\ell^{2,f}$.
\end{fed}
\begin{ex}
\emph{As an example of a fractal scale Hilbert space let $S^1=\R/\Z$
be the circle and
$$\mathcal{H}=\big\{W^{k,2}(S^1,\R)\big\}_{k \in \N_0}$$
the scale Hilbert space of Sobolev maps from the circle to the
reals. Here we understand that
$$W^{0,2}(S^1,\R)=L^2(S^1,\R)$$
the space of square integrable maps. An orthogonal basis for
$L^2(S^1,\R)$ is the Fourier basis $\{e_\nu\}_{\nu \in \N}$ defined
by $$e_1(t)=1,\,\,e_{2m}(t)=\sqrt{2}\sin(2\pi m t), \,\,
e_{2m+1}(t)=\sqrt{2}\cos(2\pi m t), \quad m \in \N,\,\,t \in S^1.$$
The Fourier basis is a common orthogonal basis for all
$W^{k,2}(S^1,\R)$. Indeed, one computes for $\nu,\nu' \in \N$
$$\langle e_\nu,e_{\nu'}\rangle_k=\delta_{\nu,\nu'}
\Bigg(\sum_{j=0}^k\bigg(2\pi \bigg\lfloor
\frac{\nu}{2}\bigg\rfloor\bigg)^{2j}\Bigg).$$ We conclude that
$\mathcal{H}$ is scale isomorphic to $\ell^{2,\sigma}$ where $\sigma
\colon \N \to (0,\infty)$ is the function
$$\sigma(\nu)=\nu^2+1, \quad \nu \in \N.$$}
\end{ex}
\begin{rem}
If $\mathcal{H}=\{H_k\}_{k \in \N_0}$ is a scale Hilbert space then
following H.\,Hofer, K.\,Wysocki, and E.\,Zehnder we denote for $m
\in \N_0$ by $\mathcal{H}^m$ the shifted scale Hilbert space
\begin{equation}\label{shift}
\mathcal{H}^m=\{H_{k+m}\}_{k \in \N_0}.
\end{equation}
If $\mathcal{H}$ is fractal then $\mathcal{H}^m$ is scale isomorphic
to $\mathcal{H}=\mathcal{H}^0$ for every $m \in \N_0$. On the other
hand this property is not sufficient to characterize fractal scale
Hilbert spaces as can be shown with the methods from
\cite{frauenfelder}.
\end{rem}
We next justify why we believe that the scale Hilbert spaces which
arise in semiinfinite dimensional Morse theories as introduced by
Floer \cite{floer} are expected to be fractal. This reasoning is
inspired by the paper \cite{robbin-salamon} of J.\,Robbin and
D.\,Salamon.

Recall that if $H_1$ is a vector space and $H_0$ is a Hilbert space,
then a linear operator $A \colon H_1 \to H_0$ is called
\emph{Fredholm} if the following three conditions hold
\begin{description}
 \item[(i)] $\mathrm{ker}(A)$ is finite dimensional.
 \item[(ii)] $\mathrm{im}(A)$ is a closed subspace of the Hilbert
 space $H_1$.
 \item[(iii)] The orthogonal complement $\mathrm{im}(A)^\perp \subset H_0$ is
 finite dimensional.
\end{description}
\begin{rem}
\emph{One usually requires that $H_1$ is itself a Hilbert space or
at least a Banach space and $A$ is continuous. However, note that
the Fredholm property is independent of the norm on $H_1$ and only
needs the scalar product on $H_0$. \emph}
\end{rem}
If $A \colon H_1 \to H_0$ is a Fredholm operator, then the integer
$$\mathrm{ind}(A)=\dim(\mathrm{ker}(A))-\dim(\mathrm{im}(A)^\perp)$$
is called the \emph{Fredholm index} of $A$. For the following
definition recall also the notation $\mathcal{H}^m$ for the shifted
scale Hilbert spaces from (\ref{shift}).
\begin{fed}
Let $\mathcal{H}$ be a scale Hilbert space and $A \colon
\mathcal{H}^1 \to \mathcal{H}^0$ be a linear scale continuous
operator. Then $A$ is called a \emph{scale Hessian operator} if the
following three axioms hold.
\begin{description}
 \item[Symmetry: ] The operator $A \colon H_1 \to H_0$ is symmetric,
 i.e. $\langle \xi, A \eta \rangle_0=\langle A \xi, \eta \rangle_0$
 for all $\xi, \eta \in H_1$.
 \item[Fredholm: ]The operator $A \colon H_1 \to H_0$ is a Fredholm
  operator of index $0$.
 \item[Regularity: ] If $\xi \in H_1$ and $A \xi \in H_n$ for
  $n \in \N_0$, then actually $\xi \in H_{n+1}$.
\end{description}
\end{fed}
\begin{rem}\emph{The motivation for introducing the notion of
scale Hessian operator is that in Floer theory the Hessian at a
critical point is supposed to be a scale Hessian operator. The
Hessian should obviously be symmetric. The requirement that the
Hessian is a Fredholm operator of index $0$ can be interpreted that
the critical point equation represents a well-posed problem which is
neither overdetermined nor underdetermined. The regularity axiom can
be thought of as to guarantee smoothness of the solutions of the
critical point equation.}
\end{rem}
\textbf{Theorem A: } \emph{Assume that for the scale Hilbert space
$\mathcal{H}$ a scale Hessian operator exists. Then $\mathcal{H}$ is
fractal.}
\\ \\
\textbf{Proof: } Theorem~A is an immediate consequence of
Corollary~\ref{mr} proved in the following section. \hfill $\square$

\section[Proof of the main result]{Proof of the main result}

Let $A$ be a scale Hessian operator on a scale Hilbert space
$\mathcal{H}$. We can also interpret the linear operator $A \colon
H_1 \to H_0$ as an unbounded operator $A \colon H_0 \to H_0$ with
dense domain $\dom(A)=H_1 \subset H_0$. The following Lemma tells us
that $A$ interpreted in this way is a self-adjoint operator.

\begin{lemma}\label{sa}
Assume that $H_0$ is a Hilbert space and $A$ is a symmetric
unbounded operator on $H_0$ with dense domain $\mathrm{dom}(A)=H_1
\subset H_0$. Suppose further that $A \colon H_1 \to H_0$ is a
Fredholm operator of index $0$. Then $A$ is selfadjoint.
\end{lemma}
\textbf{Proof: } We first note that the symmetry of $A$ implies that
$$\mathrm{ker}(A) \subset \mathrm{im}(A)^\perp.$$
Since the index of $A$ is $0$ by assumption we conclude that
$$\dim(\mathrm{ker}(A))=\dim(\mathrm{im}(A)^\perp)$$
and hence
\begin{equation}\label{fre}
\mathrm{ker}(A)=\mathrm{im}(A)^\perp.
\end{equation}
 Now choose
$$v \in \mathrm{dom}A^*.$$
This means that there exists $y \in H_0$ such that the following
holds
$$\langle y, w \rangle=\langle v, Aw \rangle, \quad \forall\,\, w
\in H_1.$$ Choose $\eta \in H_1$ and $\xi \in \mathrm{im}(A)^\perp$
satisfying
$$y=A\eta+ \xi.$$
For $q \in H_1 \cap \mathrm{im}(A)$ we compute
$$\langle v,Aq \rangle=\langle y,q \rangle=\langle A \eta,q \rangle+
\langle \xi, q \rangle=\langle \eta, A q \rangle.$$ We conclude that
$$\langle v-\eta,Aq \rangle=0, \quad \forall\,\,q
\in H_1 \cap \mathrm{im}(A)=H_1 \cap \mathrm{ker}(A)^\perp.$$ But
this is equivalent to the assertion
$$\langle v-\eta,w \rangle=0, \quad \forall\,\,w \in
\mathrm{im}(A).$$ We deduce using (\ref{fre})
$$v-\eta \in \mathrm{im}(A)^\perp=\mathrm{ker}(A).$$
Hence
$$v \in \eta +\mathrm{ker}(A) \subset \mathrm{dom}(A).$$
In particular, $A$ is selfadjoint. This finishes the proof of the
Lemma. \hfill $\square$
\\ \\
If $A$ is a scale Hessian operator on a scale Hilbert space
$\mathcal{H}$, we can endow $H_1=\dom(A)$ with another inner product
defined by
$$\langle \xi, \eta \rangle_A=\langle \xi, \eta \rangle_0+
\langle A \xi, A \eta \rangle_0, \quad \xi, \eta \in H_1.$$ The norm
$|| \cdot ||_A$ induced from $\langle \cdot, \cdot \rangle_A$ is the
graph norm on $\dom(A)$, and since $A$ is self-adjoint by
Lemma~\ref{sa} we conclude that the pair $\big(\dom(A),
\langle\cdot,\cdot \rangle_A\big)$ is complete, i.e.~itself a
Hilbert space.
\begin{lemma}\label{sc}
The two scalar products $\langle \cdot, \cdot \rangle_1$ and
$\langle \cdot, \cdot \rangle_A$ on $H_1=\dom(A)$ give rise to
equivalent norms.
\end{lemma}
\textbf{Proof: } We prove the Lemma in two Steps.
\\ \\
\textbf{Step~1: } \emph{There exists a constant $c>0$ such that
$||\xi||_1 \leq c||\xi||_A$ for all $\xi \in H_1$.}
\\ \\
We complexify the real Hilbert spaces $H_1$ and $H_0$ to complex
Hilbert spaces $H_1^\C$ and $H_0^C$. Since $A$ interpreted as
unbounded operator on $H_0^C$ is selfadjoint by Lemma~\ref{sa} we
conclude that the spectrum of $A$ is contained on the real axis. In
particular there exists $\lambda \in \C$ such that $(A-\lambda
\mathrm{id}) \colon H_1^C \to H_0^C$ is invertible. Since $A \colon
\mathcal{H}^1 \to \mathcal{H}^0$ is scale continuous we conclude
that $A-\lambda \mathrm{id}$ is a continuous bijective map from
$H^\C_1$ to $H^C_0$. Hence by the open mapping theorem its inverse
is also continuous. In particular, there exists $c_0>0$ such that
$$||(A-\lambda \mathrm{id})^{-1}\eta||_1 \leq c_0 ||\eta||_0, \quad
\forall\,\,\eta \in H_0^\C.$$ Hence we estimate for $\xi \in H_1$
\begin{eqnarray*}
||\xi||_1&=&||(A-\lambda \mathrm{id})^{-1}(A-\lambda \mathrm{id})
\xi||_1\\
&\leq&c_0||(A-\lambda \mathrm{id}) \xi||_0\\
&\leq&c_0\big(||A\xi||_0+||\lambda \xi||_0\big)\\
&\leq&\max\{c_0,|\lambda|c_0\}\big(||\xi||_0+||A\xi||_0\big).
\end{eqnarray*}
Hence Step~1 follows with $c=\max\{c_0,|\lambda|c_0\}$.
\\ \\
\textbf{Step~2: } \emph{We prove the Lemma.}
\\ \\
By Step~1 the identity map is a continuous map from
$(H_1,\langle\cdot,\cdot\rangle_A)$ to
$(H_1,\langle\cdot,\cdot\rangle_1)$. Hence by the open mapping
theorem the two scalar products are equivalent. \hfill $\square$
\begin{fed}
An unbounded selfadjoint operator $A$ on a Hilbert space $H_0$ is
called \emph{cocompact} if the inclusion $I \colon
\big(\dom(A),\langle \cdot, \cdot \rangle_A\big) \to
\big(H_0,\langle\cdot,\cdot \rangle_0\big)$ is a compact operator.
\end{fed}

\begin{lemma}\label{coco}
Assume that $A$ is a cocompact selfadjoint operator on a Hilbert
space $H_0$. Then the spectrum of $A$ is discrete and consists of
eigenvalues of finite multiplicity.
\end{lemma}
\textbf{Proof: } We set $H_1=\dom(A)$ and endow it with the scalar
product $\langle \cdot, \cdot \rangle_A$. As in the proof of
Lemma~\ref{sc} we again complexify the Hilbert spaces $H_0$ and
$H_1$ to complex Hilbert spaces $H_0^\C$ and $H_1^C$ and note that
there exists $\lambda \in \C$ such that $A_\lambda=A-\lambda
\mathrm{id} \colon H_1^\C \to H_0^\C$ is invertible. We consider the
operator
$$B_\lambda \colon H^\C_0 \to H^C_0$$
defined by
$$B_\lambda=I \circ A_\lambda^{-1}.$$
By the open mapping theorem $A_\lambda^{-1}$ is continuous and $I$
is compact by assumption. We conclude that $B_\lambda$ is a compact
operator. We next show that $B_\lambda$ is normal. In order to do
that we first determine its adjoint. For $\xi_1,\xi_2 \in H_0$ we
compute with respect to the scalar product $\langle \cdot, \cdot
\rangle$ on $H_0$ using the fact that $A$ is selfadjoint
\begin{eqnarray*}
\langle B_\lambda^*\xi_1,\xi_2 \rangle&=&\langle \xi_1,B_\lambda\xi_2 \rangle\\
&=&\langle \xi_1,I A_\lambda^{-1} \xi_2 \rangle\\
&=&\langle
A_{\bar{\lambda}}A_{\bar{\lambda}}^{-1}\xi_1,A_\lambda^{-1}\xi_2\rangle\\
&=&\langle A_{\bar{\lambda}}^{-1}\xi_1,A_\lambda
A_\lambda^{-1}\xi_2\rangle\\
&=&\langle I A_{\bar{\lambda}}^{-1}\xi_1,\xi_2\rangle.
\end{eqnarray*}
We deduce
$$B_\lambda^*=I \circ A_{\bar{\lambda}}^{-1}=B_{\bar{\lambda}}.$$
Using this formula we obtain
\begin{eqnarray*}
B^*_\lambda B_\lambda&=&I A_{\bar{\lambda}}^{-1}IA_\lambda^{-1}\\
&=&I A_{\bar{\lambda}}^{-1}A_\lambda^{-1}\\
&=&I(A_\lambda A_{\bar{\lambda}})^{-1}\\
&=&I\big(A^2-2\mathrm{Re}(\lambda)A+|\lambda|^2\mathrm{id}\big)^{-1}\\
&=&B_\lambda B_\lambda^*.
\end{eqnarray*}
In particular, $B_\lambda$ is normal. Summing up, we have shown that
$B_\lambda$ is a compact, invertible, normal operator. We apply next
the spectral theorem to $B_\lambda$ and conclude that there exists
an orthonormal basis $\{e_\nu\}_{\nu \in \N}$ of $H_0^\C$ consisting
of eigenvectors of $B_\lambda$ and a zero sequence $\{\mu_\nu\}_{\nu
\in \N}$ of corresponding eigenvalues of $B_\lambda$ such that for
$\xi \in H_0$ the formula
$$B_\lambda \xi=\sum_{\nu=1}^\infty \mu_\nu \langle \xi,e_\nu\rangle
e_\nu$$ holds. If $\nu \in \N$ we get by applying $A_\lambda$ to the
eigenvalue equation
$$B_\lambda e_\nu=\mu_\nu e_\nu$$
the equation
$$A_\lambda e_\nu=\frac{1}{\mu_\nu} e_\nu.$$
By the definition of $A_\lambda$ we conclude
\begin{equation}\label{eig}
Ae_\nu=\bigg(\frac{1}{\mu_\nu}+\lambda\bigg)e_\nu=:\gamma_\nu e_\nu.
\end{equation}
In particular, if $\xi \in H_1^\C$ we have the formula
\begin{equation}\label{sto}
A \xi=\sum_{\nu=1}^\infty \gamma_\nu \langle \xi,e_\nu\rangle e_\nu.
\end{equation}
Since $A$ is selfadjoint it follows that $\gamma_\nu$ is real for
each $\nu \in \N$, so that we can replace the basis $\{e_\nu\}_{\nu
\in \N}$ by a real orthonormal basis of $H_0$ such that (\ref{sto})
continues to hold. This finishes the proof of the Lemma. \hfill
$\square$
\\ \\
We continue the notation of the proof of the previous Lemma. After
reordering the basis vectors $\{e_\nu\}_{\nu \in \N}$ we can assume
without loss of generality that the function $\nu \mapsto
|\gamma_\nu|$ is mononone. We define $f_A \in
\widetilde{\mathcal{F}}$ by the formula
$$f_A(\nu)=1+\gamma_\nu^2, \quad \nu \in \N.$$
\begin{cor}
Assume that $A$ is a cocompact selfadjoint operator on a Hilbert
space $H_0$. Then the scale Hilbert pair $\big(H_0,\dom(A)\big)$ is
scale isometric to $(\ell^2,\ell^2_{f_A})$.
\end{cor}
\textbf{Proof: }Identify $H_0$ with $\ell^2$ using an orthonormal
basis for which (\ref{sto}) holds with monotone increasing
$|\gamma_\nu|$. We compute for $\nu,\mu \in \N$
$$\langle e_\nu,e_\mu \rangle =\langle e_\nu, e_\mu \rangle
+\langle A e_\nu, A e_\mu \rangle= \delta_{\nu\mu}+\gamma_\nu
\gamma_\mu \delta_{\nu \mu}.$$ This proves the Corollary. \hfill
$\square$
\begin{cor}
Assume that $\mathcal{H}$ is a scale Hilbert space and $A$ is a
scale Hessian operator on $\mathcal{H}$. Then the restriction of $A$
to $\mathcal{H}^2$ is a scale Hessian operator on $\mathcal{H}^1$.
\end{cor}
\textbf{Proof: }The symmetry axiom and the regularity axiom for
$A|_{\mathcal{H}^2}$ are trivial. The only nontrivial axiom to check
is the Fredholm axiom. By the regularity axiom we obtain
$$\dom(A|_{H_1})=H_2.$$
Let $\{e_\nu\}_{\nu \in \N}$ be an orthonormal basis of $H_0$ as
constructed in the proof of Lemma~\ref{coco} such that (\ref{sto})
holds true. Since the Fredholm property is further unchanged if we
replace the scalar product by an equivalent one we can thanks to
Lemma~\ref{sc} assume without loss of generality that $H_1=\dom(A)$
is endowed with the scalar product $\langle \cdot,\cdot \rangle_A$.
For $\nu \in \N$ we set
$$e_\nu^A=\frac{1}{\sqrt{f_A(\nu)}} e_\nu.$$
Then $\{e_\nu^A\}_{\nu \in \N}$ is an orthonormal basis of $H_1$. By
(\ref{sto}) we get for $\xi \in H_2$ the formula
$$A\xi=\sum_{\nu=1}^\infty \gamma_\nu \langle \xi,e_\nu^A \rangle_A
e^A_\nu.$$ In particular, if one identifies $H_1$ with $H_0$ under
the isomorphism which maps the orthonormal basis $\{e_\nu^A\}_{\nu
\in \N}$ of $H_1$ to the orthonormal basis $\{e_\nu\}_{\nu \in \N}$
of $H_0$ the operators $A|_{H_1}$ and $A$ are identified. Since $A
\colon \dom(A) \to H_0$ is a Fredholm operator of index zero, the
same has to be true for $A|_{H_1}$. This finishes the proof of the
Corollary. \hfill $\square$

\begin{cor}\label{mr}
Assume that $A$ is a scale Hessian operator on a scale Hilbert space
$\mathcal{H}$. Then $\mathcal{H}$ is scale isomorphic to
$\ell^{2,f_A}$.
\end{cor}
\textbf{Proof: } By regularity we obtain for $k \in \N_0$ unbounded
operators
$$A_k=A|_{H_k}: \dom(A_k)=H_{k+1} \to H_k.$$
In particular, we have
$$A_{k+1}=A_k|_{\dom(A_k)}.$$
By induction on Lemma~\ref{sc} we can assume that after replacing
the scale structure of $\mathcal{H}$ by a scale isomorphic scale
structure $\mathcal{H}'$, that for $\xi, \eta \in H_{k+1}$ the inner
product is given by
\begin{equation}\label{gra}
\langle \xi, \eta \rangle_{k+1}=\langle \xi, \eta\rangle_k+ \langle
A_k \xi, A_k \eta \rangle_k.
\end{equation}
From the proof of Lemma~\ref{coco} we obtain an orthogonal basis of
$H_0$ such that (\ref{sto}) holds. We can further assume without
loss of generality that the map $\nu \mapsto |\gamma_\nu|$ is
monotone. For $k \in \N_0$ we set
\begin{equation}\label{ob}
e_\nu^k=\frac{1}{f_A(\nu)^{k/2}}e_\nu.
\end{equation}
By induction on (\ref{gra}) we obtain that $\{e^k_\nu\}_{\nu \in
\N}$ is an orthonormal basis of $H_k$ and if $\xi \in
H_{k+1}=\dom(A_k)$ the formula
$$A_k \xi=\sum_{\nu=1}^\infty \gamma_\nu \langle \xi,e_\nu^k
\rangle_k e^k_\nu$$ holds true. We infer from (\ref{ob}) that
$\mathcal{H}'=\ell^{2,f_A}$. But $\mathcal{H}'$ is scale isomorphic
to $\mathcal{H}$ which finishes the proof of the Corollary. \hfill
$\square$

\end{document}